\documentclass{article}
\usepackage{amssymb}
\usepackage{amsmath}
\usepackage[english]{babel}

\newtheorem{theorem}{Theorem}

\newtheorem{corollary}[theorem]{Corollary}

\newtheorem{lemma}[theorem]{Lemma}

\newtheorem{remark}[theorem]{Remark}

\begin{document}
\title{Spectral estimates for degenerate critical levels.}
\author{Brice Camus,\\
Ruhr-Universit\"at Bochum, Fakult\"at f\"ur Mathematik,\\
Universit\"atsstr. 150, D-44780 Bochum, Germany.\\
Email : brice.camus@univ-reims.fr}
\maketitle
\begin{abstract}
\noindent We establish spectral estimates at a critical energy
level for $h$-pseudors . Via a trace formula, we compute the
contribution of isolated (non-extremum) critical points under a
condition of "real principal type". The main result holds for all
dimensions, for a singularity of any finite order and can be
invariantly expressed in term of the geometry of the singularity.
When the singularities are not integrable on the energy surface
the results are significative since the order w.r.t. $h$ of the
spectral
distributions are bigger than in the regular setting.\medskip\\
Keywords : Semi-classical analysis; Trace formula; Oscillatory
integrals.
\end{abstract}
%%%%%%%%%%%%%%%%%%%%%%%%%%%%%%%%%%%%%%%%%%%%%%%%%%%%%%%%%%%%%%%%%%%%%
%%%%%%%%%%%%%%%%%%%%%%%%%%%%%%%%%%%%%%%%%%%%%%%%%%%%%%%%%%%%%%%%%%%%%
\section{Introduction.}
If $P_{h}$ is a self-adjoint $h$-pseudo-differential operator, or
more generally $h$-admissible (see \cite{[Rob]}), acting on a
dense subset of $L^2$, a classical and accessible problem is to
study the asymptotic behavior, as $h$ tends to 0, of the spectral
distributions :
\begin{equation}
\gamma (E,\varphi,h)=\sum\limits_{|\lambda _{j}(h)-E|\leq
\varepsilon }\varphi (\frac{\lambda _{j}(h)-E}{h}),  \label{Def
trace}
\end{equation}
where the $\lambda _{j}(h)$ are the eigenvalues of $P_{h}$. Here
we suppose that the spectrum is discrete in $[E-\varepsilon
,E+\varepsilon ]$, a sufficient condition for this is given below.
One motivation is that it is in general not possible to compute
the spectrum and one has to use statistical methods to gain
spectral information. A second motivation is the existence of a
duality between spectrum of quantum operators and the classical
mechanic attached to their symbols.

If $p_0$ is the principal symbol of $P_{h}$, an energy $E$ is
regular when $\nabla p_0(x,\xi )\neq 0$ on the energy surface
$\Sigma _{E}=\{(x,\xi )\in T^{\star}\mathbb{R}^{n}\text{ }/\text{
}p_0(x,\xi )=E\}$ and critical when it is not regular. It is well
known that asymptotics of (\ref{Def trace}), as $h$ tends to 0,
are closely related to the closed trajectories of the classical
flow of $p_0$ on the surface $\Sigma_{E}$. Hence, there is a
duality between the following objects :
\begin{equation*}
\lim_{h\rightarrow 0}\gamma (E,\varphi, h)\rightleftharpoons
\{(t,x,\xi)\in\mathrm{supp}(\hat{\varphi})\times \Sigma _{E}
\text{ / } \Phi_{t}(x,\xi)=(x,\xi)\},
\end{equation*}
where $\Phi_{t}$ is the flow of the Hamiltonian vector field
$H_{p_0}=\partial _{\xi }p_0.\partial_{x} -\partial
_{x}p_0.\partial_{\xi}$. This duality has a universal character
and does not systematically requires the presence of an asymptotic
parameter, as can show the trace formulae of Selberg and
Duistermaat-Guillemin \cite{D-G}. In the semi-classical setting,
this relation was initially pointed out in the physic literature
\cite{GUT}\&\cite{BB}. From a mathematical point of view, and $E$
a regular energy, a non-exhaustive list of references concerning
this subject is \cite{BU}, \cite{PU}, \cite{P-P} and more recently
with a different approach \cite{CRR}. See also \cite{HR} for the
case of elliptic operators.

If $E$ is no more a regular value, the behavior of (\ref{Def
trace}) depends on the singularities of $p$ on $\Sigma_{E}$ which
leads to technical complications. The case of non-degenerate
critical energies, that is such that the critical-set
$\mathbb{\frak{C}}(p_0) =\{(x,\xi )\in T^{\ast
}\mathbb{R}^{n}\text{ }/\text{ }dp_0(x,\xi )=0\}$ is a compact
$C^{\infty }$ manifold with a Hessian $d^{2}p_0$ transversely
non-degenerate along this manifold has been studied first by
Brummelhuis et al. \cite{BPU}. The problem was solved for quite
general operators but for some ''small times'', i.e. for
$\mathrm{supp}(\hat{\varphi})$ contained in a neighborhood of the
origin so that the only period of the linearized flow in
$\rm{supp}(\hat{\varphi})$ was 0. Later, Khuat-Duy \cite{KhD1} has
obtained the contributions of the non-zero periods of the
linearized flow with $\rm{supp}(\hat{\varphi})$ compact, but for
Schr\"{o}dinger operators with symbol $\xi ^{2}+V(x)$ and a
non-degenerate potential $V$. Our contribution was to generalize
the result of \cite{KhD1} for more general operators but under
extra assumptions on the flow (see \cite{Cam}). Finally, in
\cite{Cam1} the case of totally degenerate extremum was treated
and the objective of this work is to study degenerate
singularities which are not an extremum of the symbol.

After a reformulation, based on the theory of Fourier integral
operators, the asymptotics of (\ref{Def trace}) can be expressed
in terms of oscillatory integrals whose phases are related to the
flow of $p_0$ on $\Sigma_E$. When $(x_0,\xi_0)$ is a critical
point of $p_0$, it is well known that the relation
$\mathrm{Ker}(d_{x,\xi}\Phi_{t}(x_0,\xi_0)-\mathrm{Id})\neq \{
0\}$ leads to the study of degenerate oscillatory integrals. In
this work we consider the case of a totally degenerate energy,
that is such that the Hessian matrix at our critical point is
zero. Hence, the linearized flow for such a critical point
satisfies $d_{x,\xi }\Phi _{t}(x_{0},\xi _{0})=\mathrm{Id}$, for
all $t\in\mathbb{R}$. A fortiori :
\begin{equation} \label{singularity}
\mathrm{Ker}(d_{x,\xi}\Phi_{t}(x_0,\xi_0)-\mathrm{Id})=T_{x_0,\xi_0}T^*\mathbb{R}^n\simeq
\mathbb{R}^{2n}, \text{ } \forall t\in\mathbb{R},
\end{equation}
and the oscillatory integrals we have to consider are totally
degenerate. In particular, it is impossible to use the stationary
phase method to determine the asymptotic behavior of Eq. (\ref{Def
trace}).

To solve this problem we will establish suitable normal forms, for
the phase functions of our oscillatory integrals, for whom it is
possible to generalize the stationary phase formula. The
construction is geometric and is independent of the dimension but
the asymptotic expansion of the related oscillatory integrals
depends on the dimension and on the order of the singularity at
the critical point. Finally, since the normal forms have a
geometrical meaning, it is possible to express invariantly the top
order coefficients of the asymptotic expansions in term of the
geometry of the singularity on the energy surface.

\section{Hypotheses and main result.}
Let $P_{h}=Op_{h}^{w}(p(x,\xi ,h))$ an $h$-pseudodifferential
operator, obtained by Weyl quantization, in the class of
$h$-admissible operators with symbol $p(x,\xi ,h) \sim \sum
h^{j}p_{j}(x,\xi )$, i.e. there exists $p_{j}\in \Sigma
_{0}^{m}(T^{\ast }\mathbb{R} ^{n})$ and $R_{N}(h)$ such that :
\begin{equation*}
P_{h}=\sum\limits_{j<N}h^{j}p_{j}^{w}(x,hD_{x})+h^{N}R_{N}(h),\text{
} \forall N\in \mathbb{N}.
\end{equation*}
Here $R_{N}(h)$ is a bounded family of operators on
$L^{2}(\mathbb{R}^{n})$ for $h\leq h_{0}$ and :
\begin{equation*}\Sigma _{0}^{m}(T^{\ast }\mathbb{R}^{n})=\{a:T^{\ast }\mathbb{R}%
^{n}\rightarrow \mathbb{C},\text{ }\sup|\partial _{x}^{\alpha
}\partial _{\xi }^{\beta }a(x,\xi )|<C_{\alpha ,\beta }m(x,\xi
),\text{ }\forall \alpha ,\beta \in \mathbb{N}^{n}\},
\end{equation*}
where  $m$ is a tempered weight on $T^{\ast }\mathbb{R}^{n}$. For
a detailed exposition on $h$-admissible operators we refer to the
book of Robert \cite{[Rob]}. In particular, $p_{0}(x,\xi )$ is the
principal symbol of $P_{h}$ and $p_{1}(x,\xi )$ the sub-principal
symbol. We note $\Phi _{t}:=\textrm{exp} (tH_{p_0})$, the
Hamiltonian flow of $H_{p_0}=\partial _{\xi }p_0
.\partial_{x}-\partial _{x}p_0 .\partial_{\xi }$ and
$\Sigma_E=p_0^{-1}(E)\subset T^{*}\mathbb{R}^n$ the energy
surfaces of $p_0$.

We study here asymptotics of the spectral distributions :
\begin{equation}
\gamma (E_{c},\varphi,h)=\sum\limits_{\lambda _{j}(h)\in \lbrack
E_{c}-\varepsilon ,E_{c}+\varepsilon ]}\varphi (\frac{\lambda
_{j}(h)-E_{c}}{h}), \text{ }h\rightarrow 0^+,\label{Objet trace}
\end{equation}
under the hypotheses $(H_{1})$ to $(H_{4})$ given below. We use
here the notation $E_c$ to recall that this energy will be chosen
critical.
\medskip\\
$(H_{1})$\textit{ The symbol of }$P_{h}$ \textit{is real and there
exists } $\varepsilon_{0}>0$ \textit{ such that the set }%
$p_0^{-1}([E_c-\varepsilon_{0},E_c+\varepsilon_{0}])$\textit{\ is
compact in $T^*\mathbb{R}^n$.}
\begin{remark}\label{spectrum} \rm{By Theorem 3.13 of \cite{[Rob]} the spectrum $\sigma (P_{h})\cap
[E_{c}-\varepsilon ,E_{c}+\varepsilon ]$ is discrete and consists
in a sequence $\lambda _{1}(h)\leq \lambda _{2}(h)\leq ...\leq
\lambda _{j}(h)$ of eigenvalues of finite multiplicities, if
$\varepsilon$ and $h$ are small enough. A fortiori, $(H_1)$
insures that $\Sigma_{E_c}$ is compact.}
\end{remark}
To simplify notations we write $z=(x,\xi)$ for any point of the
phase space.\medskip\\
$(H_{2})$\textit{ On }$\Sigma _{E_c}$\textit{, }$p_0$\textit{ has
a unique critical point }$z_{0}=(x_{0},\xi _{0})$\textit{ and near
}$z_{0}$ :
\begin{equation*}
p_0(z)=E_{c}+\sum\limits_{j=k}^{N}\mathfrak{p}_{j}(z)+\mathcal{O}(||(z-z_{0})||^{N+1}),\text{
} k>2,
\end{equation*}
\textit{where the functions }$\mathfrak{p}_{j}$\textit{\ are homogeneous of degree }$j$%
\textit{\ in }$z-z_{0}.$
\begin{remark}\rm{Strictly speaking, one could consider $k=2$ with
$\mathrm{supp}(\hat{\varphi})$ small. But there is nothing new
here since this case is precisely treated in \cite{BPU}.}
\end{remark}
$(H_{3})$\textit{ We have }$\hat{\varphi}\in C_{0}^{\infty
}(\mathbb{R}).$\ \medskip\\
Since we are interested in the contribution of the fixed point
$z_{0}$, to understand the new phenomenon it suffices to study a
local problem  :
\begin{equation}
\gamma _{z_{0}}(E_{c},\varphi,h)=\frac{1}{2\pi }\mathrm{Tr}\int\limits_{\mathbb{R}}e^{i%
\frac{tE_{c}}{h}}\hat{\varphi}(t)\psi ^{w}(x,hD_{x})\mathrm{exp}(-\frac{i}{h}%
tP_{h})\Theta (P_{h})dt.
\end{equation}
Here $\Theta $ is a function of localization near the critical
energy surface $\Sigma_{E_c}$ and $\psi \in C_{0}^{\infty
}(T^{\ast }\mathbb{R}^{n})$ has an appropriate support near
$z_{0}$. Rigorous justifications are given in section 3 for the
introduction of $\Theta (P_{h})$ and in section 4 for $\psi
^{w}(x,hD_{x})$. The case of a critical point which is not an
extremum is quite difficult, in particular because the singularity
is transferred on the blow up of the critical point. To obtain a
reasonable problem we consider the following hypothesis inspired
by H\"ormander's
real principal condition for distributions : \medskip\\
$(H_{4})$ \textit{We have }$\nabla \mathfrak{p}_{k}\neq 0$\textit{
on the set }$C (\mathfrak{p}_k)=\{ \theta \in
\mathbb{S}^{2n-1}\text{ / } \mathfrak{p}_{k}(\theta )=0
\}$\textit{.}
\begin{remark} \rm{Contrary to the case of a local extremum \cite{Cam1},
$z_0$ is not isolated on $\Sigma_{E_c}$. This imposes to study the
classical dynamic in a neighborhood of $z_0$.}
\end{remark}
By homogeneity we have $\nabla \mathfrak{p}_{k}\neq 0$ on the cone
$\{z\in T^*\mathbb{R}^n \backslash \{0\} \text{ /
}\mathfrak{p}_k(z)=0\}$, but mainly $(H_4)$ will be used on the
unit sphere. $C_{\mathfrak{p}_k}$ is a smooth manifold
of codimension 1 which can be equipped with an invariant (Liouville) measure.\medskip\\
\textbf{Construction of a geometrical measure.} To state the main
result clearly, we explain how to construct a relevant measure on
the line. By $(H_4)$ we can locally construct on
$C(\mathfrak{p}_k)$ the $2n-2$ dimensional Liouville form $dL$ via
:
\begin{equation}
dL(\theta)\wedge d \mathfrak{p}_k (\theta) =d\theta, \text{ }
\forall \theta \in C(\mathfrak{p}_k).
\end{equation}
Here the differential are w.r.t. coordinates on
$\mathbb{S}^{2n-1}$, $d\theta$ is the standard surface density and
we note again $\mathfrak{p}_k (\theta)$ for the restriction of
$\mathfrak{p}_k$ on $\mathbb{S}^{2n-1}$. The form $dL$ induces a
density on $C(\mathfrak{p}_k)$ and by continuity we can extend the
construction to close surfaces
$\mathfrak{p}_k(\theta)=\varepsilon$ for $\varepsilon$ small
enough. Accordingly we define the integrated density as :
\begin{equation}
\mathrm{Lvol}(u)=\int\limits_{\{\mathfrak{p}_k(\theta)=u \}}
|dL(\theta)|.
\end{equation}
An alternative definition is to compute the volume in
$\mathbb{S}^{2n-1}$ of the pullback :
\begin{equation}
\mathfrak{p}_k^{-1}([0,x])=\{\theta \in \mathbb{S}^{n-1}\text{ / }
\mathfrak{p}_k(\theta)\in [0,x] \},
\end{equation}
and to interpret the result as a measure :
\begin{equation}
\mathrm{V}(\mathfrak{p}_k^{-1}([0,x]))=\int\limits_{0}^{x}
\mathrm{Lvol}(s)ds.
\end{equation}
This relation is known in geometry as the co-area formula. The
most important point, that will be exploited in this work is that
$(H_4)$ insures that $\mathrm{Lvol}(u)$ is smooth near the origin.
In other words, viewing $\mathrm{Lvol}(u)$ as a distribution we
obtain that $0\notin
\mathrm{singsupp}(\mathrm{Lvol}(u))$.\medskip\\
With these objects the new contributions to the trace formula are
given by :
\begin{theorem}
\label{Main1} Under hypotheses $(H_{1})$ to $(H_{4})$ we obtain
the existence of a full asymptotic expansion :
\begin{equation}
\gamma_{z_0}(E_c,\varphi,h)\sim
h^{\frac{2n}{k}-n}\sum\limits_{m=0,1} \sum\limits_{j=0}^\infty
h^{\frac{j}{k}}\log(h)^m \Lambda_{j,m}(\varphi),
\end{equation}
where the logarithms only occur when $(2n+j)/k\in \mathbb{N}^*$
and $\Lambda_{j,m}\in \mathcal{S}'(\mathbb{R})$.\\
As concerns the leading term we obtain :\\
(1) If $k>2n$ (non-integrable singularity on $\Sigma_{E_c}$) we
have :
\begin{equation*}
\gamma _{z_{0}}(E_{c},\varphi,h)\sim h^{\frac{2n}{k}-n}\Lambda
_{0,0}(\varphi )+\mathcal{O}(h
^{\frac{2n+1}{k}-n}\mathrm{log}(h)),\text{ as }h\rightarrow 0,
\end{equation*}
where the distributional coefficient $\Lambda _{0,0}(\varphi )$ is
given by :
\begin{equation}\label{resultat1}
\frac{1}{(2\pi)^{n}k} ( \left\langle
t_{+}^{\frac{2n}{k}-1},\varphi(t)\right\rangle \int\limits_{\{
\mathfrak{p}_k \geq 0\}} |\mathfrak{p}_k (\theta)
|^{-\frac{2n}{k}}d\theta + \left\langle
t_{-}^{\frac{2n}{k}-1},\varphi(t)\right\rangle \int\limits_{\{
\mathfrak{p}_k \leq 0\}} |\mathfrak{p}_k (\theta)
|^{-\frac{2n}{k}}d\theta ).
\end{equation}
(2) If the ratio $2n/k=q$ is an integer we obtain logarithmic
contributions :
\begin{equation*}
\gamma _{z_{0}}(E_{c},\varphi,h)\sim
h^{\frac{2n}{k}-n}\mathrm{log}(h)\Lambda_{0,1} (\varphi
)+\mathcal{O}(h ^{\frac{2n}{k}-n}),\text{ as }h\rightarrow 0,
\end{equation*}
where :
\begin{equation} \Lambda_{0,1} (\varphi )
=\frac{1}{(2\pi )^n}\frac{d^{q-1}\mathrm{LVol}}{du^{q-1}}
(0)\int\limits_{\mathbb{R}} |t|^{q-1} \varphi(t) .
\end{equation}
(3) For $2n>k$ and $2n/k\notin \mathbb{N}$ the asymptotic is as in
1) with the modified distributions :
\begin{equation*}
\left\langle t_{+}^{\frac{2n}{k}-1}, \varphi(t)\right\rangle
\left\langle \frac{\tilde{d}^{2n}}{\tilde{d}u^{2n}}
u^{2n-\frac{2n}{k}}_{+}, \mathrm{Lvol} \right\rangle +\left\langle
t_{-}^{\frac{2n}{k}-1}, \varphi(t)\right\rangle \left\langle
\frac{\tilde{d}^{2n}}{\tilde{d}u^{2n}} u^{2n-\frac{2n}{k}}_{-},
\mathrm{Lvol} \right\rangle,
\end{equation*}
where the derivatives w.r.t. $w$ are normalized distributional
derivatives.
\end{theorem}
The meaning of normalized derivative it that one choose the
normalization :
\begin{equation}
\frac{\tilde{d}^{2n}}{\tilde{d}w^{2n}}=\prod\limits_{j=1}^{2n}\frac{1}{j-\frac{2n}{k}}
\frac{d^{2n}}{dw^{2n}},
\end{equation}
so that the distributional derivatives satisfy :
\begin{equation}\label{normalized derivatives}
\left\langle \frac{\tilde{d}^{2n}}{\tilde{d}w^{2n}}
w^{2n-\frac{2n}{k}}_{\pm}, f(w)\right\rangle = \left\langle
w^{-\frac{2n}{k}}_{\pm}, f(w)\right\rangle,
\end{equation}
for all $f\in C_0^\infty$ with $f=0$ in a neighborhood of the
origin. With the method we employ here this normalization appears
naturally in the expansion. The distributional bracket involving
$\mathrm{Lvol}$ is detailed in section 6.
\begin{remark}
\rm{ In cases $\textit{(1)}$\&$\textit{(3)}$ the remainders
$\mathcal{O}(h ^{\frac{2n+1}{k}-n}\mathrm{log}(h))$ can be
replaced by $\mathcal{O}(h ^{\frac{2n+1}{k}-n})$ under the only
condition that $(2n+1)/k\notin \mathbb{N}$. An interesting point
is that the degree $k$ of the singularity on $T^*\mathbb{R}^n$ has
the inverse scaling property on the spherical blow-up since the
new degree is $-2n/k$.}
\end{remark}
Results $\textit{(3)}$\&$\textit{(2)}$ for $q\geq 2$ are not
intuitive and are certainly difficult to be reached without
geometry. In particular one has to work in the dual since both
Fourier transforms w.r.t. $t$ are distributional. In
$\textit{(3)}$, the order $2n$ is arbitrary and the result is the
same for any derivative of order greater than $\mathrm{E}(2n/k)$.
Viewing the top order coefficients of the trace as distributions,
i.e. :
\begin{equation}
\gamma_{z_0}(E_c,\varphi,h)\sim f(h) \left\langle \gamma,
\varphi\right\rangle, \text{ }h\rightarrow 0,
\end{equation}
in all cases at hand we obtain :
\begin{corollary} \label{singular support}
Under the previous assumptions, $\mathrm{singsupp}(\gamma)=\{0\}$.
\end{corollary}
A similar result presumably holds for all terms of the expansion
since the asymptotic involves distributions $|t|^\alpha
\log(|t|)$. Note that Corollary \ref{singular support} is not
obvious in view of Eq. (\ref{singularity}). Also, it must be
pointed out that results $\textit{(1)}$\&$\textit{(2)}$ for $q=1$
are bigger than the standard estimate for non-critical energies
for which one obtain :
\begin{equation} \label{regular}
\gamma(E,\varphi,h)\sim \frac{h^{1-n}}{(2\pi)^n} \hat{\varphi}(0)
\mathrm{Vol}(\Sigma_E),
\end{equation}
where $\mathrm{Vol}(\Sigma_E)$ is the usual Liouville volume of
the regular (compact) surface of energy $E$. Hence the presence of
non-integrable singularities on the energy surface has a
significative spectral effect which can perhaps be measured by
eigenfunctions estimates as in \cite{BPU}. On the other side, for
an integrable singularity, we obtain the global result :
\begin{corollary} Under the conditions of Theorem \ref{Main1}, if
$k<2n$ we have :
\begin{equation}
\gamma(E_c,\varphi,h) \sim \frac{h^{1-n}}{(2\pi)^n}
\hat{\varphi}(0) \mathrm{Vol}(\Sigma_{E_c}),\text{ as }
h\rightarrow 0.
\end{equation}
\end{corollary}
This result is a consequence of Theorem \ref{Main1} and of the
results of \cite{BPU} section 4 to which we refer for a detailed
proof. Finally, an interesting problem that we have not
investigated here is the problematic of small $h$ and $|E-E_c|$
estimates for a non-integrable singularity on $\Sigma_{E_c}$ : in
this setting $\mathrm{Vol}(\Sigma_E)$ of Eq.(\ref{regular})
diverges as $E\rightarrow E_c$. See \cite{BPU} for more details.
\section{Oscillatory representation.}
The construction below is more or less classical. We recall here
important facts and results concerning the approximation of the
propagator by Fourier integral operators, or FIO, depending on a
parameter $h$. We recall that :
\begin{equation*} \gamma (E_{c},\varphi,h)
=\sum\limits_{\lambda _{j}(h)\in I_{\varepsilon }}\varphi
(\frac{\lambda _{j}(h)-E_{c}}{h}), \text{ } I_{\varepsilon }
=[E_{c}-\varepsilon ,E_{c}+\varepsilon ],
\end{equation*}
with $\hat{\varphi}\in C_{0}^{\infty }(\mathbb{R})$ and
$p_{0}^{-1}(I_{\varepsilon_{0}})$ compact in $T^{\ast
}\mathbb{R}^{n}$. In this setting the spectrum of $P_{h}$ is
discrete in $I_{\varepsilon }$ for $h>0$ small enough and
$\varepsilon<\varepsilon_{0}$ (see Remark \ref{spectrum}) and the
sum is well defined. We localize near $E_{c}$ with a cut-off
$\Theta \in C_{0}^{\infty }(]E_{c}-\varepsilon ,E_{c}+\varepsilon
\lbrack )$, $\Theta =1$ near $E_{c}$ and $0\leq \Theta \leq 1$ on
$\mathbb{R}$. The associated decomposition is :
\begin{equation*}
\gamma (E_{c},\varphi,h)=\gamma _{1}(E_{c},\varphi,h)+\gamma
_{2}(E_{c},\varphi,h),
\end{equation*}
with :
\begin{gather}\label{gamma_1 O(h infiny)}
\gamma _{1}(E_{c},\varphi,h)=\sum\limits_{\lambda _{j}(h)\in
I_{\varepsilon }}(1-\Theta )(\lambda _{j}(h))\varphi
(\frac{\lambda _{j}(h)-E_{c}}{h}),\\
\gamma _{2}(E_{c},\varphi,h)=\sum\limits_{\lambda _{j}(h)\in
I_{\varepsilon }}\Theta (\lambda _{j}(h))\varphi (\frac{\lambda
_{j}(h)-E_{c}}{h}).
\end{gather}
A classical result, see e.g. \cite{Cam1}, is that the sum of Eq.
(\ref{gamma_1 O(h infiny)}) satisfies :
\begin{lemma}
$\gamma _{1}(E_{c},\varphi,h)=\mathcal{O}(h^{\infty })$, as
$h\rightarrow 0$.\label{S1(h)=Tr}
\end{lemma}
Consequently, asymptotic behaviors of $\gamma (E_{c},\varphi,h)$
and $\gamma _{2}(E_{c},\varphi,h)$ are equivalent modulo
$\mathcal{O}(h^{\infty })$. By inversion of the Fourier transform
we obtain :
\begin{equation*}
\Theta (P_{h})\varphi (\frac{P_{h}-E_{c}}{h})=\frac{1}{2\pi}\int\limits_{%
\mathbb{R}}e^{i\frac{tE_{c}}{h}}\hat{\varphi}(t)\mathrm{exp}(-\frac{i}{h}%
tP_{h})\Theta (P_{h})dt.
\end{equation*}
Since the trace of the left hand-side is exactly $\gamma
_{2}(E_{c},\varphi,h)$, we have :
\begin{equation}
\gamma _{2}(E_{c},\varphi,h)
=\frac{1}{2\pi}\mathrm{Tr}\int\limits_{\mathbb{R}}e^{i\frac{tE_{c}}{h}}
\hat{\varphi}(t)\mathrm{exp}(-\frac{i}{h}tP_{h})\Theta (P_{h})dt.
\label{Trace S2(h)}
\end{equation}
Observe that Eq. (\ref{Trace S2(h)}) generalizes the Poisson
summation formula. \medskip\\
Let be $U_{h}(t)=\mathrm{exp}(-\frac{it}{h}P_{h})$, the evolution
operator. For each integer $N$ we can approximate $U_{h}(t)\Theta
(P_{h})$, modulo $\mathcal{O}(h^{N})$, by a Fourier
integral-operator, or FIO, depending on a parameter $h$. Let
$\Lambda$ be the Lagrangian manifold associated to the flow of
$p_0$, i.e. :
\begin{equation*}
\Lambda =\{(t,\tau ,x,\xi ,y,\eta )\in T^{\ast }\mathbb{R}\times T^{\ast }%
\mathbb{R}^{n}\times T^{\ast }\mathbb{R}^{n}:\tau =p_0(x,\xi
),\text{ }(x,\xi )=\Phi _{t}(y,\eta )\}.
\end{equation*}
\begin{theorem}
$U_{h}(t)\Theta (P_{h})$ is an $h$-FIO associated to $\Lambda$.
There exist $U_{\Theta ,h}^{(N)}(t)$ with integral kernel in
H\"ormander's class $I(\mathbb{R}^{2n+1},\Lambda )$ and
$R_{h}^{(N)}(t)$ bounded, with a $L^{2}$-norm uniformly bounded
for $0<h\leq 1$ and $t$ in a compact subset of $\mathbb{R}$, such
that
$U_{h}(t)\Theta(P_{h})=U_{\Theta,h}^{(N)}(t)+h^{N}R_{h}^{(N)}(t)$.
\end{theorem}
We refer to Duistermaat \cite{DUI1} for a proof of this theorem.
\begin{remark}
\rm{By a theorem of Helffer\&Robert (Theorem 3.11 and Remark 3.14
of \cite{[Rob]}), $\Theta (P_{h})$ is an $h$-admissible operator
with a symbol supported in $p_0^{-1}(I_{\varepsilon }).$ This
allows us to consider only oscillatory-integrals with compact
support.}
\end{remark}
For the control of the remainder, associated to $R_{h}^{(N)}(t)$,
we use :
\begin{corollary}
Let $\Theta _{1}\in C_{0}^{\infty }(\mathbb{R})$ such that $\Theta
_{1}=1$ on $\rm{supp}(\Theta )$ and $\rm{supp}(\Theta _{1})\subset
I_{\varepsilon }$, then $\forall N\in \mathbb{N}$ :
\begin{equation*}
\mathrm{Tr}(\Theta (P_{h})\varphi (\frac{P_{h}-E_{c}}{h}))=\frac{1}{2\pi }%
\mathrm{Tr}\int\limits_{\mathbb{R}}\hat{\varphi}(t)e^{\frac{i}{h}%
tE_{c}}U_{\Theta ,h}^{(N)}(t)\Theta
_{1}(P_{h})dt+\mathcal{O}(h^{N}).
\end{equation*}
\end{corollary}
For a proof of this result, based on the cyclicity of the trace, see \cite{Cam1} or \cite{[Rob]}.\medskip\\
If $(x_{0},\xi _{0})\in \Lambda $ and $\varphi (x,\theta )\in
C^{\infty }(\mathbb{R}^{k}\times \mathbb{R}^{N})$ parameterizes
$\Lambda $ in a sufficiently small neighborhood $U$ of $(x_{0},\xi
_{0})$, then for each $u_{h}\in I(\mathbb{R}^{k},\Lambda )$ and
$\chi \in C_{0}^{\infty }(T^{\ast }\mathbb{R}^{k})$,
$\rm{supp}(\chi )\subset U,$ there exists a
sequence of amplitudes $a_{j}=a_{j}(x,\theta )\in C_{0}^{\infty }(\mathbb{R}%
^{k}\times \mathbb{R}^{N})$ such that for all $L\in\mathbb{N}$ :
\begin{equation}
\chi ^{w}(x,hD_{x})u_{h}=\sum\limits_{-d\leq j<L}h^{j}I(a_{j}e^{\frac{i}{h}%
\varphi })+\mathcal{O}(h^{L}).
\end{equation}
We will use this remark with the following result of H\"{o}rmander
(see \cite {HOR1}, tome 4, proposition 25.3.3). Let be $(T,\tau
,x_{0},\xi _{0},y_{0},-\eta _{0})\in \Lambda _{\mathrm{flow}}$,
$\eta _{0}\neq 0$, then near this point there exists, after
perhaps a change of local coordinates in $y$ near $y_{0},$ a
function $S(t,x,\eta )$ such that :
\begin{equation}
\phi(t,x,y,\eta )=S(t,x,\eta )-\left\langle y,\eta \right\rangle ,
\end{equation}
parameterizes $\Lambda _{\mathrm{flow}}$. In particular this
implies that :
\begin{equation*}
\{(t,\partial _{t}S(t,x,\eta ),x,\partial _{x}S(t,x,\eta
),\partial _{\eta }S(t,x,\eta ),-\eta )\}\subset
\Lambda_{\mathrm{flow}} ,
\end{equation*}
and that the function $S$ is a generating function of the flow,
i.e. :
\begin{equation}
\Phi _{t}(\partial _{\eta }S(t,x,\eta ),\eta ) =(x,\partial
_{x}S(t,x,\eta )). \label{Gene}
\end{equation}
Moreover, $S$ satisfies the Hamilton-Jacobi equation :
\begin{equation*}
\left\{
\begin{array}{c}
\partial _{t}S(t,x,\eta )+ p_0(x,\partial
_{x}S(t,x,\eta ))=0, \\
S(0,x,\xi)=\left\langle x,\xi \right\rangle .
\end{array}
\right.
\end{equation*}
Mainly, we will apply this result locally near $(x_{0},\xi
_{0})=(y_{0},\eta _{0})$, our unique fixed point of the flow on
the energy surface $\Sigma _{E_{c}}$. If $\xi _{0}=0$ we can
replace the operator
$P_{h}$ by $e^{\frac{i}{h}\left\langle x,\xi _{1}\right\rangle }P_{h}e^{-\frac{i}{h}%
\left\langle x,\xi _{1}\right\rangle }$ with $\xi _{1}\neq 0.$ The
spectrum is the same since the new operator has symbol $p(x,\xi
-\xi _{1},h)$ with critical point $(x_{0},\xi _{1})$, $\xi
_{1}\neq 0$.

Hence, for each $N\in \mathbb{N}^{*}$ and modulo an error
$\mathcal{O}(h^{N-d})$, the localized trace $\gamma
_{2}(E_{c},\varphi,h)$ of Eq. (\ref{Trace S2(h)}) can be written
as :
\begin{equation}
\gamma _{2}(E_{c},\varphi,h)=\sum\limits_{j<N}(2\pi h)^{-d+j}\int\limits_{\mathbb{%
R\times R}^{2n}}e^{\frac{i}{h}(S(t,x,\xi )-\left\langle x,\xi
\right\rangle +tE_{c})}a_{j}(t,x,\xi )\hat{\varphi}(t)dtdxd\xi .
\label{gamma1 OIF}
\end{equation}
To obtain the right power $-d$ of $h$ we apply results of
Duistermaat \cite{DUI1} (following here H\"{o}rmander for the FIO,
see \cite {HOR2} tome 4) concerning the order. An
$h$-pseudo-differential operator obtained by Weyl quantization :
\begin{equation*}
(2\pi h)^{-\frac{N}{2}}\int\limits_{\mathbb{R}^{N}}a(\frac{x+y}{2},\xi )e^{%
\frac{i}{h}\left\langle x-y,\xi \right\rangle }d\xi ,
\end{equation*}
is of order 0 w.r.t. $1/h$. Since the order of $U_{h}(t)\Theta (P_{h})$ is $-%
\frac{1}{4}$, we have :
\begin{equation}
\psi ^{w}(x,hD_{x})U_{h}(t)\Theta (P_{h})\sim
\sum\limits_{j<N}(2\pi
h)^{-n+j}\int\limits_{\mathbb{R}^{n}}a_{j}(t,x,y,\eta )e^{\frac{i}{h}%
(S(t,x,\eta )-\left\langle y,\eta \right\rangle )}dy.
\label{operateur d'evolution}
\end{equation}
Multiplying Eq. (\ref{operateur d'evolution}) by $\hat{\varphi}%
(t)e^{\frac{i}{h}tE_{c}}$ and passing to the trace we find Eq.
(\ref{gamma1 OIF}) with $d=n$ and where we write again
$a_{j}(t,x,\eta )$ for $a_{j}(t,x,x,\eta )$.

To each element $u_{h}$ of $I(\mathbb{R}^{k},\Lambda )$ we can
associate a principal symbol $e^{\frac{i}{h}S}\sigma
_{\mathrm{princ}}(u_{h})$, where $S$ is a function on $\Lambda $
such that $\xi dx=dS$ on $\Lambda .$ In fact, if
$u_{h}=I(ae^{\frac{i}{h}\varphi })$ then we have $S=S_{\varphi
}=\varphi \circ i_{\varphi }^{-1}$ and $\sigma
_{\mathrm{princ}}(u_{h})$ is a section of $|\Lambda
|^{\frac{1}{2}} \otimes M(\Lambda )$, where $M(\Lambda )$ is the
Maslov vector-bundle of $\Lambda $ and $|\Lambda |^{\frac{1}{2}}$
the bundle of half-densities on $\Lambda$. If $p_{1}$ is the
sub-principal symbol of $P_{h}$, in the global coordinates
$(t,y,\eta )$ on $\Lambda _{\mathrm{flow}}$ the half-density of
the propagator $U_{h}(t)$ is given by :
\begin{equation}
\exp (i\int\limits_{0}^{t}p_{1}(\Phi _{s}(y,-\eta ))ds)|dtdyd\eta |^{\frac{1%
}{2}}.\label{demi densite}
\end{equation}
This expression is related to the resolution of the first
transport equation for the propagator, for a proof we refer to
Duistermaat and H\"{o}rmander \cite{D-H}.
\section{Classical dynamic near the equilibrium.}
The function $S$ of Eq. (\ref{gamma1 OIF}) is related to the
classical mechanic and we obtain a link with the dynamic generated
by $p_0$. Precisely, a critical point of the oscillatory integral
of Eq. (\ref{gamma1 OIF}) satisfies the equations :
\begin{equation*}
\left\{
\begin{matrix}
E_{c}=-\partial _{t}S(t,x,\xi ), \\
x=\partial _{\xi }S(t,x,\xi ), \\
\xi =\partial _{x}S(t,x,\xi ),
\end{matrix}
\right. \Leftrightarrow \left\{
\begin{array}{c}
p_0(x,\xi )=E_{c}, \\
\Phi _{t}(x,\xi )=(x,\xi ),
\end{array}
\right.
\end{equation*}
where the right hand side defines a closed trajectory of the flow
inside $\Sigma _{E_{c}}$. Accordingly, the non-stationary phase
method shows that $\gamma _{2}(E_{c},\varphi,h)$ is asymptotically
determined by the closed orbits of the flow on $\Sigma _{E_{c}}$.
We recall that we are mainly interested in the contribution of the
equilibrium $z_0$. With $\psi \in C_{0}^{\infty }(T^{\ast
}\mathbb{R}^{n})$, $\psi =1\text{ near }z_{0}$, we write :
\begin{gather*}
\gamma _{2}(E_{c},\varphi,h) =\frac{1}{2\pi }\mathrm{Tr}\int\limits_{\mathbb{R}}e^{i%
\frac{tE_{c}}{h}}\hat{\varphi}(t)\psi ^{w}(x,hD_{x})\mathrm{exp}(-\frac{i}{h}%
tP_{h})\Theta (P_{h})dt\\
+\frac{1}{2\pi }\mathrm{Tr}\int\limits_{\mathbb{R}}e^{i\frac{tE_{c}}{h}}\hat{%
\varphi}(t)(1-\psi
^{w}(x,hD_{x}))\mathrm{exp}(-\frac{i}{h}tP_{h})\Theta (P_{h})dt.
\end{gather*}
With the additional hypothesis of having a clean flow, the second
term can be fully treated by the semi-classical trace formula on a
regular level. We also observe that the contribution of the first
term is micro-local and allows to introduce local coordinates near
$z_{0}$. To distinguish the contribution of $z_{0}$ from eventual
closed trajectories we use the following result on the dynamic.
\begin{lemma}\label{dynamic}
For all $T>0$ there exists a neighborhood $U_{T}$ of the critical
point such that $\Phi _{t}(z)\neq z$ for all $z\in U_{T}\backslash
\{z_{0}\}$ and for all $t\in ]\text{-}T,0[\cup ]0,T[$.
\end{lemma}
For a proof of this Lemma, based mainly on a result of Yorke
\cite{Yor}, see \cite{Cam1}.\medskip\\
With $\mathrm{supp}(\hat{\varphi})$ compact we can choose $\psi$
such that Lemma \ref{dynamic} holds on $\mathrm{supp}(\psi )$ for
all $t\in \rm{supp}(\hat{\varphi})$. Hence, on the support of
$\psi $ there is two contributions :\\
1) Points $(t,x,\xi )=(0,x,\xi )$ for $(x,\xi )\in \Sigma
_{E_{c}}$.\\
2) Points $(t,x,\xi )=(t,z_{0})$ for $t\in
\rm{supp}(\hat{\varphi})$.\\
Now, we restrict our attention to the second contribution. In the
following, until further notice, the derivatives $d$ are
derivatives w.r.t. initial conditions $z$ (resp. $x,\xi$). If
$z_0$ is an equilibrium of the vector field $X$ then the
linearized flow $d\exp(tX)(z_0)$ is the flow of the linearized
vector field $dX(z_0)$. In our setting, the linearized operator is
zero and we have :
\begin{equation} \label{linear flow}
d\Phi _{t}(z_{0})=\mathrm{exp}(0)=\textrm{Id},\text{ }\forall t.
\end{equation}
Moreover, with $(H_{2})$ we clearly obtain that :
\begin{equation}d^{j}\Phi
_{t}(z_{0})=0,\text{ }\forall t,\text{ }\forall j\in \{2,..,k-2\}.
\end{equation}
The non-zero terms of the Taylor expansion of the flow are
computed by :
\begin{lemma}
\label{TheoFormule de récurence du flot}Let be $z_{0}$ an
equilibrium of the $C^{\infty}$ vector field $X$ and $\Phi _{t}$
the flow of $X$. Then for all $m\in \mathbb{N}^{\ast }$, there
exists a polynomial map $P_{m}$, vector valued and of degree at
most $m$, such that :
\begin{equation*}
d^{m}\Phi _{t}(z_{0})(z^{m})=d\Phi
_{t}(z_{0})\int\limits_{0}^{t}d\Phi _{-s}(z_{0})P_{m}(d\Phi
_{s}(z_{0})(z),...,d^{m-1}\Phi _{s}(z_{0})(z^{m-1}))ds.
\end{equation*}
\end{lemma}
See \cite{Cam} or \cite{Cam1} for a proof. In our setting, the
$(k-1)$-jet of $p_0$ is flat in $z_0$ and we have
$P_{k-1}(y_1,...,y_{k-2})=d^{k-1}H_{p_0}(z_{0})(y_1^{k-1})$. Here,
for any vector $u$ the notation $u^l$ stands for $(u,...,u)$
repeated $l$-times, the same convention is used below. In view of
Eq. (\ref{linear flow}) by integration form 0 to $t$ we obtain :
\begin{equation}
d^{k-1}\Phi_{t}(z_{0})(z^{k-1})=\int\limits_{0}^{t}
d^{k-1}H_{p_0}(z_{0})(z^{k-1})ds=td^{k-1}H_{\mathfrak{p}_{k}}(z_{0})(z^{k-1}).
\label{derive ordre k-1 du flot}
\end{equation}
This provides an explicit formula for the germ of $\Phi_t$ in
$z_0$ :
\begin{equation}\label{germ flow}
\Phi_t(z)=z
+\frac{1}{(k-1)!}d^{k-1}\Phi_{t}(z_{0})(z^{k-1})+\mathcal{O}(||z||^k).
\end{equation}
We describe now more precisely the singularities of our phase
function. Without loss of generality, we can assume that $z_0=0$
\begin{lemma}\label{structure phase} Near the origin we
have :
\begin{equation}
S(t,x,\xi )-\left\langle x,\xi \right\rangle
+tE_{c}=-t(\mathfrak{p}_{k}(x,\xi )+R_{k+1}(x,\xi
)+tG_{k+1}(t,x,\xi )), \label{forme phase}
\end{equation}
where $R_{k+1}(x,\xi ) =\mathcal{O}(||(x,\xi )||^{k+1})$ and
$G_{k+1}(t,x,\xi)=\mathcal{O}(||(x,\xi )||^{k+1})$, uniformly for
$t$ in a compact subset of $\mathbb{R}$.
\end{lemma}
\textit{Proof.} In view of Eq. (\ref{germ flow}), we search our
local generating function as :
\begin{equation*}
S(t,x,\xi)=-tE_c+\left\langle x,\xi \right\rangle+S_{k}(t,x,\xi)+
\mathcal{O}(||(x,\xi )||^{k+1}),
\end{equation*}
where $S_k$ is homogeneous of degree $k$ w.r.t. $(x,\xi)$. Let $J$
be the matrix of the usual symplectic form. Comparing terms of
degree $k-1$ in the implicit relation
$\Phi_t(\partial_{\xi}S(t,x,\xi),\xi)=(x,\partial_{x}S(t,x,\xi))$
provides :
\begin{equation*}
J\nabla S_k(t,x,\xi)=\frac{1}{(k-1)!}d^{k-1}\Phi
_{t}(0)((x,\xi)^{k-1}).
\end{equation*}
By homogeneity and with Eq. (\ref{derive ordre k-1 du flot}) we
obtain :
\begin{equation*}
S_k(t,x,\xi)=\frac{1}{k!}\left\langle
(x,\xi),tJd^{k-1}H_{\mathfrak{p}_k}(x,\xi)^{k-1}
\right\rangle=-t\mathfrak{p}_{k}(x,\xi).
\end{equation*}
As concerns the remainders, we have $S(0,x,\xi)= \left\langle
x,\xi \right\rangle $, so that :
\begin{equation*}
S(t,x,\xi)-\left\langle x,\xi \right\rangle=tF(t,x,\xi),
\end{equation*}
where $F$ is smooth in a neighborhood of $(x,\xi)=0$. Now, the
Hamilton-Jacobi equation imposes that $F(0,x,\xi)=-p_0(x,\xi)$ and
we have :
\begin{equation*}
R_{k+1}(x,\xi)=p_0(x,\xi)-E_c-\mathfrak{p}_k(x,\xi)=\mathcal{O}(||(x,\xi)||^{k+1}).
\end{equation*}
Finally, the time dependant remainder can be written :
\begin{equation*}
S(t,x,\xi)-S(0,x,\xi)-t\partial_t S(0,x,\xi)=\mathcal{O}(t^2),
\end{equation*}
since by construction this term is of order
$\mathcal{O}(||(x,\xi)||^{k+1})$ we get the desired result when
$t$ is in a compact subset of $\mathbb{R}$. $\hfill{\blacksquare}$
\begin{remark}\rm{With Lemma \ref{TheoFormule de récurence du flot},
one can compute explicitly terms of higher degree for $S$ and
$\Phi_t$. But we do not need them for the present contribution
because of some considerations of homogeneity below.}
\end{remark}
\section{Normal forms and oscillatory integrals.}
Retaining only the coefficient of highest degree w.r.t. $h$ in Eq.
(\ref{gamma1 OIF}), we have to study the asymptotic behavior of
oscillatory integrals :
\begin{equation} \label{IO}
I(\frac{1}{h})=\int\limits_{\mathbb{R}\times T^*\mathbb{R}^n}
e^{\frac{i}{h}(S(t,x,\xi)-\left\langle x,\xi\right\rangle-tE_c)}
a(t,x,\xi)dtdxd\xi, \text{ } h\rightarrow 0^+.
\end{equation}
Here, we have temporary discarded the factor $(2\pi h)^{-n}$ to
avoid it's constant repetition in the calculations. Since the
contribution we study is local, we can work with local coordinates
and we identify locally our neighborhood of the critical point in
$T^{\ast }\mathbb{R}^{n}$ with an open of $\mathbb{R}^{2n}$. With
$z=(x,\xi)\in \mathbb{R}^{2n}$, we define :
\begin{equation}
\Psi(t,z)=\Psi(t,x,\xi)=S(t,x,\xi )-\left\langle x,\xi
\right\rangle+tE_{c}. \label{defphase}
\end{equation}
The next Lemma provides a resolution of singularities for $\Psi$
w.r.t. $C(\mathfrak{p}_k)$.
\begin{lemma}\label{FN1}
Assume $P_h$ satisfies conditions $(H_{2})$ and $(H_{4})$. For all
$t$ in a compact, after a blow-up w.r.t. $z$ in a neighborhood of
$z_0$, there exists local coordinates $\eta$ such that :
\begin{gather*}
\Psi(t,z) \simeq -\eta_{0}\eta_{1}^{k},\text{ in all
directions where }\mathfrak{p}_{k}>0,\\
\Psi(t,z) \simeq + \eta_{0}\eta_{1}^{k},\text{ in all
directions where }\mathfrak{p}_{k}<0,\\
\Psi(t,z) \simeq -\eta _{0}\eta_{1}^{k}\eta_{2},\text{locally near
} C(\mathfrak{p}_{k}).
\end{gather*}
\end{lemma}
\noindent\textit{Proof.} We can assume that $z_{0}=0$. To perform
the blow-up, we use polar coordinates $z=(r,\theta )$, $\theta \in
\mathbb{S}^{2n-1}(\mathbb{R})$. By Lemma \ref{structure phase},
near $z_0$ we have :
\begin{gather*}
\Psi(t,z) \simeq -tr^{k}(\mathfrak{p}_{k}(\theta
)+rR_{k+1}(r,\theta )+trG_{k+1}(t,r,\theta )),\\
\simeq -tr^{k}(\mathfrak{p}_{k}(\theta )+Q(t,r,\theta)),
\end{gather*}
where $\mathfrak{p}_{k}(\theta )$ is the restriction of
$\mathfrak{p}_{k}$ on $\mathbb{S}^{2n-1}$ and $Q(t,0,\theta)=0$.
If $\mathfrak{p}_{k}(\theta_{0})\neq 0$ and $t$ in a compact, then
for $r<r_0$ we have $\mathfrak{p}_{k}(\theta )+Q(t,r,\theta)\neq
0$. We define :
\begin{gather*}
(\eta_{0},\eta_{2},...,\eta_{2n})(t,r,\theta) =(t,\theta _{1},...,\theta _{2n-1}), \\
\eta_{1}(t,r,\theta) = r|\mathfrak{p}_{k}(\theta
)+Q(t,r,\theta)|^{\frac{1}{k}}.
\end{gather*}
In these coordinates the phase becomes $-\eta_{0}\eta_{1}^{k}$ if
$\mathfrak{p}_{k}(\theta_0)$ is positive (resp.
$\eta_{0}\eta_{1}^{k}$ for a negative value). Near $\theta_0$, we
have :
\begin{equation*}
\frac{\partial \eta_{1}}{\partial r} (t,0,\theta) =
|\mathfrak{p}_{k}(\theta )|^{\frac{1}{k}}\neq 0, \text{ } \forall
t,
\end{equation*}
hence, the corresponding Jacobian satisfies $|J\eta|(t,0,\theta
)=|\mathfrak{p}_{k}(\theta )|^{\frac{1}{k}}\neq 0$.

Now, let $\theta_0 \in C(\mathfrak{p}_{k})$. Up to a permutation,
we can suppose that
$\partial_{\theta_{1}}\mathfrak{p}_{k}(\theta_{0})\neq 0$. We
accordingly choose the new coordinates :
\begin{gather*}
(\eta_{0},\eta_{1},\eta_{3},...,\eta_{2n})(t,r,\theta)=(t,r,\theta
_{2},...,\theta _{2n-1}), \\
\eta_{2}(t,r,\theta)=\mathfrak{p}_{k}(\theta)+Q(t,r,\theta),
\end{gather*}
which are locally admissible since $|J\eta
|(t,0,\theta_0)=|\partial_{\theta_{1}}\mathfrak{p}_{k}(\theta_{0})|\neq
0$. Finally, lemma follows by compactness of $C(\mathfrak{p}_{k})$
. \hfill{$\blacksquare$}\medskip\\
If necessary, we can shrink the support of $\psi$ to obtain the
existence of the normal forms inside
$\mathrm{supp}(\hat{\varphi})\times \mathrm{supp}(\psi)$. We
define now a partition of unity associated to $\mathfrak{p}_k$. We
pick a family of functions $\phi_{j}\in
C^\infty(\mathbb{S}^{2n-1})$ such that :
\begin{equation*}
C(\mathfrak{p}_{k})\subset \bigcup\limits_{j}
\textrm{supp}(\phi_{j}),\text{ } \sum\limits_{j} \phi_{j}=1 \text{
near } C(\mathfrak{p}_{k}).
\end{equation*}
We can also choose each $\mathrm{supp}(\phi_j)$ small enough so
that normal forms of Lemma \ref{FN1} exist in
$[0,r_0]\times\textrm{supp}(\phi_{j})$. Clearly, this family can
be chosen finite and we obtain a partition of unity with
$\phi_{0}=1-\sum \phi_{j}$. The support of $\phi_{0}$ is not
connected and we define $\phi_{0}^{+}$, with
$\mathfrak{p}_k(\theta)>0$ on $\mathrm{supp}(\phi_{0}^{+})$,
similarly we define $\phi_{0}^{-}$ where $\mathfrak{p}_k<0$, so
that $\phi_{0}=\phi_{0}^{+}+\phi_{0}^{-}$.
\begin{remark}\rm{The family $(\phi_j)$ depends only on
$\mathfrak{p}_k$, e.g. we can impose $\sum \phi_{j}=0$ for
$\mathfrak{|p}_k(\theta)|\geq \varepsilon>0$. This point is useful
for the globalization in section 6.}
\end{remark}
Let be $\lambda=h^{-1}$. We accordingly split up the integral of
Eq. (\ref{IO}) to obtain :
\begin{gather*}
I_{\pm}(\lambda)=\int\limits_{\mathbb{R\times R}_{+}\times
\mathbb{S}^{2n-1}}e^{i\lambda \Psi (t,r,\theta
)}\phi_{0}^{\pm}(\theta) a(t,r\theta )r^{2n-1}dtdrd\theta\\
=\int\limits_{\mathbb{R\times R}_{+}} e^{-i\lambda(\pm \eta
_{0}\eta _{1}^{k})}A_{0}^{\pm}(\eta _{0},\eta _{1})d\eta _{0}d\eta
_{1}=\int\limits_{\mathbb{R}_{+}}\hat{A}_{0}^{\pm}(\pm \lambda
\eta _{1}^{k},\eta _{1}),
\end{gather*}
respectively for the directions where $\mathfrak{p}_k(\theta)>0$
and $\mathfrak{p}_k(\theta)<0$. Here the notation $\hat{A}$ stands
for the Fourier transform w.r.t. the first argument and is also
used below. Similarly, the covering of $C(\mathfrak{p}_k)$ gives :
\begin{gather*}
I_{j}(\lambda)=\int\limits_{\mathbb{R\times R}_{+}\times
\mathbb{S}^{2n-1}}e^{i \lambda \Psi (t,r,\theta
)}\phi_{j}(\theta) a(t,r\theta )r^{2n-1}dtdrd\theta \\
=\int\limits_{\mathbb{R}\times \mathbb{R}_{+}\times\mathbb{R}}
e^{-i \lambda \eta _{0}\eta _{1}^{k}\eta_{2}}A_{j}(\eta _{0},\eta
_{1},\eta_{2})d\eta _{0}d\eta _{1}d\eta_{2}=\int\limits_{
\mathbb{R}_{+}\times\mathbb{R}} \hat{A}_{j}(\eta
_{1}^{k}\eta_{2},\eta _{1},\eta_{2})d\eta _{1}d\eta_{2}.
\end{gather*}
These new amplitudes are obtained by pullback and integration,
i.e. :
\begin{gather}
A_{0}^{\pm}(\eta _{0},\eta _{1})=\int \eta ^{\ast
}(\phi_{0}^{\pm}(\theta)a(t,r\theta
)r^{2n-1}|J\eta |)d\eta _{2}...d\eta _{2n}, \label{ampli1}\\
A_{j}(\eta _{0},\eta _{1},\eta_{2})=\int \eta ^{\ast
}(\phi_{j}(\theta)a(t,r\theta )r^{2n-1}|J\eta |)d\eta _{3}...d\eta
_{2n}. \label{ampli2}
\end{gather}
With $C(\mathfrak{p}_k)$ compact, our oscillatory integral can be
written as a finite sum :
\begin{equation}
I(\lambda)=I_{+}(\lambda)+I_{-}(\lambda) +\sum\limits_{j=0}^{L}
I_{j}(\lambda), \text{ } \lambda=h^{-1},
\end{equation}
where each term of the r.h.s. will be treated by elementary
methods.
\begin{remark}\label{degres amplitude}
\rm{By pullback of the measure $r^{2n-1}dr$, we have
$A_{0}^{\pm}=\mathcal{O}(\eta _{1}^{2n-1})$ near $\eta _{1}=0$.
Same remark for $A_{j}=\mathcal{O}(\eta _{1}^{2n-1})$ near $\eta
_{1}=0$. This point is important, since Lemmas \ref{Theo IO 1ere
carte}, \ref{Theo IO 2eme carte} below involve Dirac-distributions
w.r.t. $\eta_1$.}
\end{remark}
\textbf{Expansion of the related oscillatory integrals.}\\
We end this section with two results on asymptotics. In fact, to
save a lot of computations we will take benefit of the linear term
$\eta_0$ in our normals forms. This approach is more economic than
the strategy proposed in \cite{WON}. But the reader must keep in
mind that the method of \cite{WON} can be applied in a more
general setting. The next elementary Lemma can be found in
\cite{Cam1} and allows to expand both integrals
$I_{\pm}(\lambda)$.
\begin{lemma} \label{Theo IO
1ere carte} For any $a\in C_{0}^{\infty
}(\mathbb{R}\times\mathbb{R}_{+})$ we have :
\begin{equation}
\int\limits_{0}^{\infty }\hat{a}(\lambda \eta _{1}^{k},\eta
_{1})d\eta _{1}\sim \sum\limits_{j=0}^{\infty }\lambda
^{-\frac{j+1}{k}}c_{j}(a), \text{ }\lambda\rightarrow +\infty,
\end{equation}
where the distributional coefficients are :
\begin{equation*}
c_{j}=\frac{1}{k}\frac{1}{j!}(\mathcal{F}(x_{+}^{\frac{j+1-k}{k}})(\eta_0)
\otimes\delta _{0}^{(j)}(\eta_1)),\text{ }x_{+}=\max (x,0).
\end{equation*}
\end{lemma}
\begin{remark}
\rm{The same result holds for a phase $-\eta _{0}\eta _{1}^{k}$ if
we replace terms $x_{+}$ by $x_{-}$ in Lemma \ref{Theo IO 1ere
carte}. For oscillatory integrals on $\mathbb{R}^2$ one can
conclude by splitting the domain of integration and 2 applications
of Lemma \ref{Theo IO 1ere carte}.}
\end{remark}
For $a\in C_0^{\infty}(\mathbb{R}\times
\mathbb{R}_+\times\mathbb{R} )$, we define the family of
elementary fiber integrals :
\begin{gather}
I_{n,k}(\lambda)=\int\limits_{0}^{\infty
}(\int\limits_{\mathbb{R}}a(\lambda y_1^k
y_2,y_1,y_2)dy_2)y_1^{2n-1}dy_1.
\end{gather}
\begin{lemma}\label{Theo IO 2eme carte}
There exists a sequence of distributions $(D_{j,p})$ such that :
\begin{equation}\label{DA singular}
I_{n,k}(\lambda)\sim \sum\limits_{p=0,1}
\sum\limits_{j\in\mathbb{N},\text{ } j\geq 2n} D_{j,p}(a) \lambda
^{-\frac{j}{k}} \log(\lambda)^p, \text{ as } \lambda\rightarrow
\infty,
\end{equation}
where the logarithms only occur when $(j/k)$ is an integer.\\
As concerns the leading term, if $(2n/k)\notin \mathbb{N}^{*}$ we
obtain :
\begin{equation} \label{equivalent non-integer}
I_{n,k}(\lambda)= \lambda^{-\frac{2n}{k}}d(a)+\mathcal{O}(\lambda
^{-\frac{2n+1}{k}}\mathrm{log}(\lambda)),
\end{equation}
with :
\begin{gather*}
d(a)=C_{n,k} \int\limits_{0}^{\infty}\int\limits_{0}^{\infty}
t^{\frac{2n}{k}-1} y_2^{2n-\frac{2n}{k}} \left(
\partial^{2n}_{y_2} a(t,0,y_2)+\partial^{2n}_{y_2}
a(-t,0,-y_2)\right) dy_2dt
\end{gather*}
But when $2n/k=q\in\mathbb{N}^{*}$, we have :
\begin{equation} \label{equivalent integer}
I_{n,k}(\lambda)= \lambda^{-\frac{2n}{k}}\log(\lambda)
\frac{1}{k}\int\limits_{\mathbb{R}} |t|^{q-1}
\partial_{y_2}^{q-1} a(t,0,0)dt
+\mathcal{O}(\lambda^{-\frac{2n}{k}}).
\end{equation}
\end{lemma}
\begin{remark}\label{remark on remainders}\rm{The remainder of Eq. (\ref{equivalent non-integer}) can be optimized
to $\mathcal{O}(\lambda ^{-\frac{2n+1}{k}})$ if $(2n+1)/k$ is not
an integer and is optimal otherwise.}
\end{remark}
\noindent\textit{Proof.} By a standard density argument we can
assume that the amplitude is of the form
$a(s,y_1,y_2)=f(s)b(y_1,y_2)$. The justification is that our
coefficients below are computed by continuous linear functionals,
i.e. distributions. We define the Melin transforms of $f$ as :
\begin{equation}
M_{\pm}(z)= \int\limits_{0}^{\infty} s^{z-1} f(\pm s)ds.
\end{equation}
We split-up $I_{n,k}$ as $J_{+}$ and $J_{-}$ by separating
integrations $y_2>0$ and $y_2<0$. By Melin inversion formula, we
accordingly obtain :
\begin{equation}\label{analytic1}
J_{+}(\lambda)= \frac{1}{2i\pi }\int\limits_{\gamma} M_{+}(z)
\lambda^{-z} \int\limits_{\mathbb{R}_{+}^2} (y_1 y_2^{k})^{-z}
b(y_1,y_2) y_1^{2n-1} dy_1dy_2 dz,
\end{equation}
where $\gamma=c+i\mathbb{R}$ and $0<c<k^{-1}$. Similarly we have :
\begin{equation}\label{analytic2}
J_{-}(\lambda)= \frac{1}{2i\pi }\int\limits_{\gamma} M_{-}(z)
\lambda^{-z} \int\limits_{\mathbb{R}_{+}^2} (y_1 y_2^{k})^{-z}
b(y_1,-y_2) y_1^{2n-1} dy_1dy_2 dz.
\end{equation}
The existence of a full asymptotic expansion is a direct
consequence of :
\begin{lemma}\label{poles extensions}
The family of distributions $z\mapsto (y_1 y_2^{k})^{-z}$ on
$C_0^{\infty}(\mathbb{R}_{+}^2)$ initially defined in the domain
$\Re(z)< k^{-1}$ is meromorphic on $\mathbb{C}$ with poles :
$z_{j,k} = j/k$, $j\in\mathbb{N}^{*}$. These poles are of order 2
when $z_{j,k}\in \mathbb{N}^*$ and of order 1 otherwise.
\end{lemma}
\textit{Proof.} We form the Bernstein-Sato polynomial $b_k$
attached to our problem :
\begin{gather*}
T(y_2y_1^{k})^{1-z}:=\frac{\partial}{\partial y_2}
\frac{\partial^{k}}{\partial y_1^{k}}
(y_2y_1^{k})^{1-z}=b_k(z)(ty_2y_1^{k})^{-z}, \\
b_k(z)=(1-z)\prod\limits_{j=1}^{k}(j-kz).
\end{gather*}
If $\Re(z)<k^{-1}$, $(k+1)$-integrations by parts yield :
\begin{equation*}
\int\limits_{\mathbb{R}_{+}^2} (y_1 y_2^{k})^{-z}
f(y_1,y_2)dy_1dy_2=\frac{(-1)^{k+1}}{b_k(z)}
\int\limits_{\mathbb{R}_{+}^2} (y_1 y_2^{k})^{1-z}
(Tf)(y_1,y_2)dy_1dy_2.
\end{equation*}
Now the integral in the r.h.s. is analytic in $\Re (z)<1+k^{-1}$.
After $p$ iterations the poles, with their orders, can be read off
the rational functions :
\begin{equation}
\mathfrak{R}_p(z) =\prod\limits_{l=1}^p \frac{1}{b_k(z-l)}.
\end{equation}
This gives the result since $p$ can be chosen arbitrary large.
$\hfill{\blacksquare}$\medskip\\
Hence, the following functions are meromorphic on $\mathbb{C}$ :
\begin{equation}
\mathfrak{g}^{\pm}(z)=\int\limits_{\mathbb{R}_{+}^2} (y_1
y_2^{k})^{-z} b(y_1,\pm y_2)dy_1dy_2.
\end{equation}
A classical result, see e.g. \cite{Bl-Hand}, is that
$M_{\pm}(c+ix)\in \mathcal{S}(\mathbb{R}_x)$ when $c \notin
-\mathbb{N}$. If we shift the path of integration $\gamma$ to the
right in our integral representation, Cauchy's residue method
provides the asymptotic expansion. For any $d>c$, outside of the
poles, we obtain that :
\begin{equation*}
\int\limits_{c+i\mathbb{R}} \lambda^{-z} M_{\pm}(z) \mathfrak{g}^{\pm}(z) dz-%
\int\limits_{d+i\mathbb{R}} \lambda^{-z} M_{\pm}(z) \mathfrak{g}^{\pm}(z) dz\\%
=\sum\limits_{c<z_{j,k}<d}
\mathrm{res}(\lambda^{-z}M_{\pm}\mathfrak{g}^{\pm})(z_{j,k}).
\end{equation*}
Since $d$ is not a pole the second integral can be estimated via :
\begin{equation}
|\int\limits_{d+i\mathbb{R}} \lambda^{-z} M_{\pm}(z)
\mathfrak{g}^{\pm}(z) dz|\leq C(f,b)
\lambda^{-d}=\mathcal{O}(\lambda^{-d}),
\end{equation}
where, for each $d$, the constant $C$ involves the $L^1$-norm of a
finite number derivatives of $b$. This will indeed lead to an
asymptotic expansion with precise remainders. Applying this method
to $J_{+}(\lambda)$ and $J_{-}(\lambda)$ we obtain the existence
of a full asymptotic expansion of the form :
\begin{equation}
I_{n,k}(\lambda) \sim \sum\limits_{p=0,1}
\sum\limits_{j\in\mathbb{N}*} C_{j,p} \lambda^{-\frac{j}{k}}
\log(\lambda)^p.
\end{equation}
Moreover, by Lemma \ref{poles extensions}, these logarithms only occur when $j/k$ is integer.\medskip\\
\textbf{Computation of the leading term.}\\
To avoid unnecessary discussions and calculations below, we remark
that we can commute the polynomial weight of Eqs.
(\ref{analytic1},\ref{analytic2}) via :
\begin{gather}
T ((y_1 y_2)^{1-z} y_1^{2n-1})=\mathfrak{b}(z)(y_2y_1^{k})^{-z}y_1^{2n-1}, \\
\mathfrak{b}(z)=(1-z)\prod\limits_{j=1}^{k}(j-kz+2n-1).
\end{gather}
By iteration, we obtain that the poles are the rational numbers :
\begin{equation*}
z_{p,j,k,n}=p+\frac{j+2n-1}{k},\text{ } j\in[1,...,k], \text{ }
p\in\mathbb{N}.
\end{equation*}
By examination of the integrals w.r.t. $y_1$, no residue
contributes before :
\begin{equation}
z_{0}= \frac{2n}{k}.
\end{equation}
This result explains the effect of the dimension $n$ and justifies
fully Eq. (\ref{DA singular}). To reach $z_0$ we need
$\mathrm{E}(2n/k)+1$ iterations. But, by analytic continuation,
any bigger integer is acceptable. For the computation of residuum
below we use :
\begin{gather}
\lambda^{-z}M_{+}(z) \mathfrak{B}_n (z)
\int\limits_{\mathbb{R}_{+}^{2}} (y_1^k y_2)^{2n-z} y_1^{2n-1}
T^{2n} b(y_1,y_2) dy_1dy_2,\\
\mathfrak{B}_n(z)=\prod\limits_{l=0}^{2n-1}
\frac{1}{\mathfrak{b}(z-l)}.
\end{gather}
This choice of $2n$ iterations is arbitrary but avoids a lot of calculations.\\
\textit{a) Case of $z_0$ simple pole.}\\
Here the residue can be computed by the limit $z\rightarrow z_0$.
We find :
\begin{equation*}
\lambda^{-\frac{2n}{k}} \lim\limits_{z\rightarrow \frac{2n}{k}}
(z-\frac{2n}{k})\mathfrak{B}_n(z) M_{+}(\frac{2n}{k})
\int\limits_{\mathbb{R}_{+}^{2}} (y_1^k y_2)^{2n-\frac{2n}{k}}
y_1^{2n-1} T^{2n} b(y_1,y_2) dy_1dy_2.
\end{equation*}
In particular, we can compute the integral w.r.t. $y_1$ via :
\begin{equation*}
\int\limits_{0}^{\infty} y_1^{2kn-1} \partial^{2kn}_{y_1}
(\partial^{2n}_{y_2} b(y_1,y_2)) dy_1 =(2kn-1)!\partial^{2n}_{y_2}
b(0,y_2).
\end{equation*}
A similar result holds for $J_{-}(\lambda)$ and we obtain :
\begin{gather}
J_{+}(\lambda)=\lambda^{-\frac{2n}{k}} C_{n,k}
M_{+}(\frac{2n}{k})
\int\limits_{0}^{\infty} y_2^{2n-\frac{2n}{k}} (\partial^{2n}_{y_2} b)(0,y_2)dy_2+R_1(\lambda),\\
J_{-}(\lambda)=\lambda^{-\frac{2n}{k}} C_{n,k} M_{-}(\frac{2n}{k})
\int\limits_{0}^{\infty} y_2^{2n-\frac{2n}{k}}
(\partial^{2n}_{y_2} b)(0,-y_2)dy_2+R_2(\lambda).
\end{gather}
Here $C_{n,k}$ is the canonical constant :
\begin{equation}\label{canonical regularization}
C_{n,k}=\frac{1}{k} \prod\limits_{j=1}^{2n}
\frac{1}{j-\frac{2n}{k}}.
\end{equation}
Finally, according to Lemma \ref{poles extensions}, each remainder
$R_j$ is of order $\mathcal{O}(\lambda^{-\frac{2n+1}{k}})$ if
$(2n+1)/k \notin \mathbb{N}$
and $\mathcal{O}(\lambda^{-\frac{2n+1}{k}}\log(\lambda))$ otherwise. This justifies Remark \ref{remark on remainders}.\medskip\\
\textit{b) Case of $z_0$ double pole.}\\
If $h$ is meromorphic with a pole of order 2 in $\xi_0$ we have :
\begin{equation*}
\mathrm{res} (h)(z_0)=\frac{1}{2} \lim\limits_{z\rightarrow z_0}
\frac{\partial}{\partial z} (z-z_0)^2 h(z).
\end{equation*}
Applying Leibnitz's rule to $\partial_z (\lambda^{-z} M_{+}(z)
\mathfrak{g}^{+}(z))$, we obtain :
\begin{equation}
J_{+}(\lambda)=B \lambda^{-\frac{2n}{k}}\log(\lambda)
+\mathcal{O}(\lambda^{-\frac{2n}{k}}),
\end{equation}
where the distribution $B$ is computed almost as before. We find :
\begin{gather*}
B=-\frac{1}{2} D_{n,k} M_{+}(\frac{2n}{k})
\int\limits_{0}^{\infty} y_2^{2n-\frac{2n}{k}} (\partial^{2n}_{y_2} b)(0,y_2)dy_2,\\
D_{n,k}= \lim\limits_{z \rightarrow \frac{2n}{k}}
(z-\frac{2n}{k})^2 \mathfrak{B}_n(z).
\end{gather*}
But $q=2n/k$ is an integer and by integration by parts we obtain :
\begin{equation}
\int\limits_{0}^{\infty} y_2^{2n-q} (\partial^{2n}_{y_2}
b)(0,y_2)dy_2 =(-1)^{q+1}(2n-q)!  \partial_{y_2}^{q-1} b(0,0).
\end{equation}
Since a similar result holds for $J_{-}(\lambda)$, we obtain the
desired result by gathering all the constants and summation.
Finally, we can extend our formulas since all coefficients in the
expansion are of the form :
\begin{equation*}
\left\langle D^j ,f\otimes b\right\rangle=\left\langle
D^j_1,f\right\rangle \left\langle D^j_2,b\right\rangle, \text{
}D^j_{1,2}\in \mathcal{D}'(\mathbb{R}).
\end{equation*}
By linearity and continuity, the result holds for a general
amplitude $a$. \hfill{$\blacksquare$}
\begin{remark} \rm{The previous method allows to compute all coefficients of the expansion.
In particular the expansion also involves logarithmic
distributions (associated to double poles). We do not detail all
these terms since they have, a priori, no invariant formulation in
the trace formula.}
\end{remark}
\section{Proof of the main result.}
We have now the desired results concerning the asymptotic behavior
of the trace. Hence, to prove the main result, it remains to
express the top order coefficients of the expansion invariantly.
Taking Remark \ref{degres amplitude} into account, to avoid
unnecessary calculations we define :
\begin{gather}
A_{0}^{\pm}(t,y_1)=y_1^{2n-1}\tilde{A}_{0}^{\pm}(t,y_1),\\
A_j(t,y_1,y_2)=y_1^{2n-1}\tilde{A}_{j}(t,y_1,y_2).
\end{gather}
Note that these definitions have no effect on the Fourier
transform w.r.t. $t$.\medskip\\
\textbf{Directions where $\mathfrak{p}_k(\theta)\neq 0$.}\\
By Lemma \ref{Theo IO 2eme carte} we obtain that the first
non-zero coefficient is obtained for $l=2n-1$ (see Remark
\ref{degres amplitude}) and is given by
\begin{equation*}
\frac{1}{k}\frac{1}{(2n-1)!}\left\langle
|\eta_{0}|_{+}^{\frac{2n-k}{k}}\otimes \delta
_{0}^{(2n-1)},\hat{A}_{0}^{+}(\eta _{0},\eta _{1})\right\rangle
=\frac{1}{k}\left\langle |\eta_{0}|_{+}^{\frac{2n-k}{k}}\otimes
\delta _{0},\hat{\tilde{A}}_{0}^{+}(\eta _{0},\eta
_{1})\right\rangle.
\end{equation*}
Since by construction :
\begin{equation}
\tilde{A}_{0}^{+}(\eta
_{0},0)=\int\limits_{\mathbb{S}^{2n-1}}a(\eta
_{0},0)\phi_{0}^{+}(\theta)|\mathfrak{p}_{k}(\theta
)|^{-\frac{2n}{k}}d\theta ,
\end{equation}
we obtain that the local contribution, associated to
$\mathrm{supp}(\phi_{0}^{+})$, is :
\begin{equation*}
\frac{1}{k}\left\langle (\mathcal{F}(x_{+}^{\frac{2n-k}{k}})
(\eta_{0}),a(\eta _{0},0)\right\rangle
\int\limits_{\mathbb{S}^{2n-1}}\phi_{0}^{+}(\theta)|\mathfrak{p}_{k}(\theta
)|^{-\frac{2n}{k}}d\theta,
\end{equation*}
and a similar result holds on $\mathrm{supp}(\phi_{0}^{-})$. Now,
since $a(t,0)=\hat{\varphi}(t)$, cf. Eq. (\ref{demi densite}), the
directions where $\mathfrak{p}_k(\theta)\neq 0$ contribute as :
\begin{gather} \label{positive}
I_{+}(\lambda) \sim \frac{1}{k} \lambda^{-\frac{2n}{k}}
\left\langle |t|_{+}^{\frac{2n-k}{k}}, \varphi (t) \right\rangle
\int\limits_{\mathbb{S}^{2n-1}}\phi_{0}^{+}(\theta)|\mathfrak{p}_{k}(\theta
)|^{-\frac{2n}{k}}d\theta,\\
\label{negative}%
I_{-}(\lambda) \sim \frac{1}{k} \lambda^{-\frac{2n}{k}}
\left\langle |t|_{-}^{\frac{2n-k}{k}}, \varphi (t) \right\rangle
\int\limits_{\mathbb{S}^{2n-1}}\phi_{0}^{-}(\theta)|\mathfrak{p}_{k}(\theta
)|^{-\frac{2n}{k}}d\theta.
\end{gather}
\textbf{Microlocal contribution of $C(\mathfrak{p}_k)$.}\newline%
According to the analysis above,
we will distinguish out the case $k$ divides $2n$.\medskip\\
$\textit{(1)}$\textit{ Case of $k>2n$, integrable singularity on the blow-up.}\\
Here $2n/k \in ]0,1[$. According to Lemma \ref{Theo IO 2eme
carte}, the contribution of $I_{j}(\lambda)$ is :
\begin{equation*}
\frac{1}{k}\lambda^{-\frac{2n}{k}}\int\limits_{\mathbb{R}_{+}^2}
|t|^{\frac{2n}{k}-1} |y_2|^{-\frac{2n}{k}}
\left(\hat{\tilde{A}}_{j}(t,0,y_2)+\hat{\tilde{A}}_{j}(-t,0,-y_2)\right)dt
dy_2 +R(\lambda).
\end{equation*}
Reminding that $y_2(t,0,\theta)=\mathfrak{p}_k(\theta)$, we obtain
:
\begin{equation*}
\int\limits_{\mathbb{R}_{+}} |y_2|^{-\frac{2n}{k}}
\tilde{A}_{j}(t,0,y_2)dy_2=a(t,0)\int\limits_{\{\mathfrak{p}_k(\theta)\geq0\}}|\mathfrak{p}_k
(\theta) |^{-\frac{2n}{k}}\phi_j(\theta) d\theta.
\end{equation*}
Via Eq. (\ref{demi densite}), by summation $I(\lambda)$ is
asymptotically equivalent to :
\begin{equation*}
\frac{\lambda^{-\frac{2n}{k}}}{k}  \left ( \left\langle
t_{+}^{\frac{2n}{k}-1},\varphi(t) \right\rangle \int\limits_{\{
\mathfrak{p}_k \geq 0\}} |\mathfrak{p}_k (\theta)
|^{-\frac{2n}{k}}d\theta + \left\langle
t_{-}^{\frac{2n}{k}-1},\varphi(t)\right\rangle \int\limits_{\{
\mathfrak{p}_k \leq 0\}} |\mathfrak{p}_k (\theta)
|^{-\frac{2n}{k}}d\theta \right ).
\end{equation*}
Note that none of the coefficients above are
equal unless $\varphi$ or $\mathfrak{p}_k$ are symmetric.\medskip\\
$\textit{(2)}$\textit{ Case of $q=2n/k$ integer.}\\
Here the contribution of each $I_{j}(\lambda)$ is dominant since
we obtain :
\begin{equation*}
I_{j}(\lambda)= \frac{1}{k}\lambda^{-q}\log (\lambda)
\int\limits_{\mathbb{R}} |t|^{q-1} \partial^{q-1}_{y_2}
\hat{\tilde{A}}_j(t,0,0)dt+\mathcal{O}(z^{-q}).
\end{equation*}
Unless $q=1$, there is no way to take the limit directly, and the
geometric properties are still hidden in the Jacobian. To reach
the result we will use the Schwartz kernel technic. Clearly, it is
enough to evaluate our derivative and to integrate w.r.t. $t$.
With $s=(s_1,s_2)\in\mathbb{R}^2$, we write the evaluation as :
\begin{equation*}
\partial^{q-1}_{y_2} \tilde{A}_j(t,0,0)=\frac{1}{(2\pi)^2} \int e^{i\left\langle
s,(y_1,y_2)\right\rangle} (is_2)^{q-1} \tilde{A}_j(t,y_1,y_2)
dy_1dy_2 ds.
\end{equation*}
Here we have used an oscillatory Schwartz kernel for $\delta_{y_1}
\otimes \delta^{q-1}_{y_2}$. This integral representation allows
to inverse our diffeomorphism to obtain :
\begin{equation*}
\partial^{q-1}_{y_2}
\tilde{A}_j(t,0,0)=\frac{1}{(2\pi)^2} \int e^{i\left\langle
s,(r,y_2(r,\theta))\right\rangle} (i s_2)^{p-1} a(t,r\theta)
\phi_j(\theta) dr d\theta ds.
\end{equation*}
If we extend the integrand by 0 for $r<0$, the normalized integral
w.r.t. $(r,s_1)$ provides $\delta_r$. By construction
$y_2(0,\theta)=\mathfrak{p}_k(\theta)$ and we accordingly have :
\begin{equation}
\partial^{q-1}_{y_2}
\tilde{A}_j(t,0,0)=a(t,0)\frac{1}{(2\pi)}
\int\limits_{\mathbb{R}\times\mathbb{S}^{n-1}} e^{i u
\mathfrak{p}_k(\theta)} (i u)^{q-1} \phi_j(\theta) d\theta du .
\end{equation}
This Fourier integral makes sense with $\mathbb{S}^{n-1}$ compact.
We define here a local version of the integrated density :
\begin{equation}
J_j(w)=\int \limits_{\{\mathfrak{p}_k(\theta)=w\}} \phi_j(\theta)
dL_w(\theta),
\end{equation}
where $dL_w$ is the density induced by the Leray-form
$dL_{\mathfrak{p}_k}$ : $d \mathfrak{p}_k \wedge
dL_{\mathfrak{p}_k}(\theta)=d\theta$. Note that all these objects
can be constructed by mean of local coordinates under the only
condition that $\mathrm{supp}(\phi_j)$ is small enough near
$C(\mathfrak{p}_k)$. Moreover, since $\mathfrak{p}_k$ is
continuous on $\mathbb{S}^{n-1}$, each $J_j(w)$ defines a
compactly supported distribution, smooth near the origin according
to $(H_4)$. The sum over all the $\phi_j$ gives the geometric
contribution :
\begin{equation}
\frac{1}{(2\pi)} \int\limits_{\mathbb{R}^2} e^{i uw} (i u)^{p-1}
\sum\limits_j J_j(w) dwdu=\frac{d^{p-1}\mathrm{Lvol}}{d
w^{p-1}}(0).
\end{equation}
By integration w.r.t. $t$ we obtain the result in general
position.
\begin{remark}
\rm{The case $2n=k$ is directly accessible with :
\begin{equation*}
I(\lambda)=\frac{1}{k} \frac{\log(\lambda)}{\lambda}
\mathrm{LVol}(0)\int\limits_{\mathbb{R}}
\varphi(t)dt+\mathcal{O}(\lambda^{-1}),
\end{equation*}
where $\mathrm{LVol}(0)$ is usual the Liouville volume of
$C(\mathfrak{p}_k)$.}
\end{remark}
$\textit{(3)}$\textit{ $k<2n$ and simple pole, non-integrable singularity.}\\
The distributional coefficients for the positive part of each
$I_{j}(\lambda)$ are :
\begin{equation*}
\left\langle \nu_{j,+},a\right\rangle
=C_{n,k}\int\limits_{0}^{\infty} \int\limits_{0}^{\infty}
t^{\frac{2n}{k}-1} y_2^{2n-\frac{2n}{k}} (\partial^{2n}_{y_2}
\tilde{A}_j)(t,0,y_2)dy_2dt.
\end{equation*}
With the same oscillatory technic as above, for
$\mathfrak{p}_k(\theta)\geq 0$ we obtain globally :
\begin{equation*}
\left\langle T_{+},a\right\rangle= C_{n,k}\left\langle
|w|_{+}^{2n-\frac{2n}{k}} , \frac{ \partial^{2n}
\mathrm{Lvol}(w)}{\partial w^{2n}}
\right\rangle\int\limits_{0}^{\infty}|t|^{\frac{2n}{k}-1}\hat{\varphi}(t)dt.
\end{equation*}
The duality bracket is well defined since $\mathrm{Lvol}$ is
compactly supported. A similar result holds for the negative part
of each $I_{j}(\lambda)$. The value of $C_{n,k}$, see Eq.
(\ref{canonical regularization}), gives the result stated in
Theorem \ref{Main1} by normalizing the distributional derivatives.
Note that our choice is convenient since $\partial_w^{2n}$ is
symmetric.
 $\hfill{\blacksquare}$\medskip\\
Now we detail the construction of the distributional bracket of
the part \textit{(3)} of Theorem \ref{Main1}. Clearly,
$\mathrm{Lvol}$ is supported in $[\inf\limits_{\mathbb{S}^{2n-1}}
\mathfrak{p}_k, \sup\limits_{\mathbb{S}^{2n-1}} \mathfrak{p}_k ]$.
Let be $\chi\in C_0^{\infty}$, $0\leq \chi \leq 1$ on $\mathbb{R}$
chosen such that $\chi=1$ near the origin and $\chi(u)=0$ for
$|u|\geq \varepsilon$, with $\varepsilon>0$ small enough. We write
the geometric contribution as :
\begin{equation*}
\left\langle U,\mathrm{Lvol} \right\rangle= \left\langle U, \chi
\mathrm{Lvol}\right\rangle+\left\langle U, (1-\chi)\mathrm{Lvol}
\right\rangle.
\end{equation*}
Away from the origin, e.g. for $u>0$, $U$ is smooth and we obtain
directly :
\begin{equation*}
C_{n,k}\left\langle \frac{d^n}{du^{2n}} u_{+}^{2n-\frac{2n}{k}},
\chi(u)\mathrm{Lvol}(u) \right\rangle
=\frac{1}{k}\int\limits_{\{f_k(\theta)
>0\}} \chi(f_k(\theta))|f_k(\theta)|^{-\frac{2n}{k}}d\theta.
\end{equation*}
Here the value of $C_{n,k}$ from Eq. (\ref{canonical
regularization}) justifies the normalization of Eq.
(\ref{normalized derivatives}). For the singular part of
$u_{+}^{-\frac{2n}{k}}$, we use the local regularity of
$\mathrm{Lvol}(u)$ near $u=0$ and integrations by parts to
conclude. Finally, using the previous trick we obtain a
formulation which does not depend on the partition of unity by a
double covering of the sphere and $\varepsilon$ small enough.
\begin{remark}\rm{The key point here is that we can put in duality the
distributions $\partial^{2n}_u |u|^{2n-\alpha}$ and
$\mathrm{Lvol}(u)$ since their singular supports are disjoints.}
\end{remark}
\textbf{Effect of the sub-principal symbol.}\\
Until now we have only considered the case of an operator given by
quantization of a symbol $p_0$. But in presence of a sub-principal
symbol $p_1$ the construction is the same. The only important
change, see Eq. (\ref{demi densite}), is that the amplitude has to
be modified by :
\begin{equation*}
a(t,z_0)= \hat{\varphi}(t) \exp( i \int\limits_{0}^{t}
p_1(\Phi_s(z_0)) ds).
\end{equation*}
But $z_0$ in a fixed point of $\Phi_t$ and a fortiori :
\begin{equation*}
a(t,z_0)= \hat{\varphi}(t) \exp( it p_1(z_0)).
\end{equation*}
Hence if $p_1(z_0)=0$, which is the case in many practical
situations, the trace formula remains the same. Otherwise, by
Fourier inversion formula, the effect is a shift on $\varphi$ by
$p_1(z_0)$ in all integral formulae.\medskip\\
\textbf{Acknowledgment.} This work was partially supported by the
\textit{SFB/TR12}, \textit{Symmetries \& Universality in
Mesoscopic Systems} and \textit{IHP-Network, Analysis \& Quantum}
ref. HPRN-CT-2002-00277.

\end{document}